\documentclass{amsart}
\usepackage{amssymb}
\usepackage{latexsym}
\usepackage{amscd}

\newtheorem{prop}{Proposition}

\begin{document}
\renewcommand{\abstractname}{Abstract}
\begin{abstract}
This note is motivated by Y.G. Oh's conjecture that the Clifford torus $L_n$ in $\mathbb{C}P^n$ minimizes volume in its Hamiltonian deformation class. We show that there exist explicit positive constants $a_n$ depending on the dimension with $a_2=3/\pi$ such that for any Lagrangian torus $L$ in the Hamiltonian class of $L_n$ we have $vol(L) \geq a_n vol (L_n)$. The proof uses the recent work of C.H. Cho \cite{Cho} on Floer homology of the Clifford tori. A formula from integral geometry enables us to derive the estimate. We wish to point out that a general lower bound on the volume of $L$ exists from the work of C. Viterbo \cite{Vit}. Our lower bound $a_2= 3/\pi$ is the best one we know.
\end{abstract}
\title[Some estimates related..]{Some estimates related to Oh's conjecture for the Clifford tori in $\mathbb{C}P^n$}     
\author{Edward Goldstein}
\maketitle
\section{Introduction}
The Clifford torus $L_n$ in $\mathbb{C}P^n$ is given in homogeneous coordinates by \[\big( (z_1 : \ldots : z_{n+1}) \big| |z_i|=|z_j| \big)\]
The Clifford torus is the only orbit of the diagonal torus action on $\mathbb{C}P^n$ which is a minimal Lagrangian submanifold, see \cite{Gold1}. It is also the only orbit which is a monotone Lagrangian submanifold, see \cite {CG}. Y.G. Oh has studied the second variation of volume of $L_n$ with respect to Hamiltonian deformations, see \cite{Oh}. He has shown that this variation is non-negative and conjectured that $L_n$ minimizes volume in its Hamiltonian deformation class. This note constitutes an effort toward verifying this conjecture. Our main tool is the recent result of Cheol-Hyun Cho \cite{Cho} which states that if $L$ is Hamiltonian equivalent to $L_n$ and if $L$ and $L_n$ intersect transversally then the number of intersection points of $L$ and $L_n$ \[\# (L \bigcap L_n) \geq 2^n\]
We will use integral geometry to study the volume of such $L$ - see also \cite{IOS} for a similar usage of integral geometry for a product of two geodesics in $S^2 \times S^2$. Our main result is that \[vol(L) \geq a_n vol (L_n)\]
with an explicit positive constant $a_n$ and $a_2=\frac{3}{\pi}$.
\section{A formula from integral geometry}
The presentation here follows R. Howard \cite{How}. In our case the group $SU(n+1)$ acts on $\mathbb{C}P^n$ with a stabilizer $K \simeq U(n)$. Thus we view $\mathbb{C}P^{n}=SU(n+1)/K$ and the Fubini-Study metric is induced from the bi-invariant metric on $SU(n+1)$. Let $P$ and $Q$ be two Lagrangian submanifolds of $\mathbb{C}P^n$. For a point $p \in P$ and $q \in Q$ we define an angle $\sigma(p,q)$ between the tangent plane $T_pP$ and $T_qQ$ as follows: First we choose some elements $g$ and $h$ in $SU(n+1)$ which move $p$ and $q$ respectively to the same point $r \in \mathbb{C}P^n$. Now the tangent planes $g_{\ast}T_pP$ and $h_{\ast}T_qQ$ are in the same tangent space $T_r \mathbb{C}P^n$ and we can define an angle between them as follows: take an orthonormal basis $u_1\ldots u_n$ for $g_{\ast}T_pP$ and an orthonormal basis $v_1 \ldots v_n$ for $h_{\ast}T_qQ$ and define \[\sigma(g_{\ast}T_pP,h_{\ast}T_qQ)= |u_1 \wedge \ldots \wedge v_n|\]
The later quantity $\sigma(g_{\ast}T_pP,h_{\ast}T_qQ)$ depends on the choices $g$ and $h$ we made. To mend this will need to average this out by the stabilizer group $K$ of the point $r$. Thus we define:
\[\sigma(p,q)=\int_{K} \sigma(g_{\ast}T_pP, k_{\ast}h_{\ast}T_qQ) dk\]
Since $SU(n+1)$ acts transitively on the Grassmanian of Lagrangian planes in $\mathbb{C}P^n$ we conclude that this angle is a constant depending just on $n$:
\[\sigma(p,q)=c_n\]
There is a following general formula due to R. Howard \cite{How}:
\[ \int_{SU(n+1)} \#(gP \bigcap Q) dg= \int_{P \times Q} \sigma(p,q) dp dq= c_n vol(P)vol(Q)\]
Thus
\begin{equation}
\label{main}
vol(P)vol(Q)= \frac{1}{c_n} \int_{SU(n+1)} \#(gP \bigcap Q)dg
\end{equation}
The quantity of interest for us is the constant $\frac{vol(SU(n+1))}{c_n}$. We'll find it using $P=Q= \mathbb{R}P^n$.
\section{The case of $\mathbb{R}P^n$}
Let $P$ be $\mathbb{R}P^n$ and let $Q$ be Hamiltonian equivalent to $P$. It is known that if $P$ and $Q$ intersect transversally then $\#(P \bigcap Q) \geq n+1$- see \cite{Giv} and also \cite{FOOO} for a more general treatment of fixed point sets of antisymplectic involutions. On the other hand if $g$ is a unitary matrix then linear algebra shows that $\#(gP \bigcap P)=n+1$ (again assuming transversality). Thus there is a proposition due to B. Kleiner:
\begin{prop}
(Kleiner) $\mathbb{R}P^n$ minimizes volume in its Hamiltonian isotopy class
\end{prop}
For our purposes we are interested in plugging the formula \ref{main} with $P=Q=\mathbb{R}P^n$. We conclude that 
\begin{equation}
\label{sup}
\frac{vol(SU(n+1))}{c_n}=\frac{vol(\mathbb{R}P^n)^2}{n+1}
\end{equation}
Let us work out the case $n=2$. The metric on $\mathbb{C}P^2$ is the quotient of the metric on $S^5$ by $S^1$-action. We have $vol(\mathbb{R}P^2)=vol(S^2)/2=2\pi$. So 
\[\frac{vol(SU(3))}{c_2}=4\pi^2/3\]
\section{The estimate for the Clifford torus}
Let $L_n \subset \mathbb{C}P^{n}$ be the Clifford torus. There is a torus $T^{n+1} \subset \mathbb{C}^{n+1}$ given by \[T^{n+1}=\big( (z_1, \ldots, z_{n+1}) \big| |z_i|=1/\sqrt{n+1} \big) \]
We have that $L_n$ is the quotient of $T^{n+1}$ by the $S^1$ action. Thus \[vol(L_n)= vol(T^{n+1})/2\pi= (2\pi/\sqrt{n+1})^{n+1}/2\pi\]
For $n=2$ we have \[vol(L_2)= 4\pi^{2}/3 \sqrt{3}\]
Let $P$ be Hamiltonian equivalent to $L_n$. From \cite{Cho} we have that for a unitary matrix $g$: $\#(gP \bigcap P) \geq 2^n$. Thus from equations \ref{main} and \ref{sup} we conclude that \[vol(P)^2 \geq 2^n \frac{vol(SU(n+1))}{c_n}=2^n \frac{vol(\mathbb{R}P^n)^2}{n+1}\]  
Let us specialize to the case $n=2$. We have \[vol(P)^2  \geq 4 \cdot 4\pi^2/3=(4\pi)^2/3 \]
Thus \[vol(P) \geq 4\pi/\sqrt{3}= \frac{3}{\pi} vol(L_2)\]

E-mail: egold@math.stanford.edu
\end{document}